\documentstyle[12pt, leqno]{article}
\begin{document}
\title{ Adiabatic Limits and Foliations} 
\author{Kefeng Liu and Weiping Zhang}
\date{}
\maketitle

\begin{abstract} We use adiabatic limits 
to study foliated manifolds. The Bott connection naturally shows up 
as the adiabatic limit of Levi-Civita connections. 
As an application, we then construct certain 
natural elliptic operators associated to the foliation and present a
 direct geometric proof of a vanshing theorem of Connes[Co],  
which extends the Lichnerowicz vanishing theorem [L]
to foliated manifolds with spin leaves, for what we call almost Riemannian
foliations. Several new vanishing theorems are also proved by using our 
method.
\end{abstract}

{\bf Introduction}

$\ $

Let $M$ be a compact spin manifold and $g^{TM}$ be a Riemannian
 metric on $TM$, the tangent bundle of $M$. 
If the scalar curvature $k_{TM}$ of $g^{TM}$ is positive, 
then a well-known theorem of Lichnerowicz [L] states that 
$\langle \widehat{A}(TM), [M]\rangle=0$. Here $\widehat{A}(TM)$ is the Hirzebruch 
$\widehat{A}$-class of $TM$. 

This result is a direct consequence of the 
Atiyah-Singer index formula [AS] and  the by now standard 
Lichnerowicz formula, which was also observed by Singer, for 
Dirac operators on spin manifolds.

Now let $(M,F)$ be a compact foliated manifold. We make the following 
assumtions: (1) $F$ is spin; (2) there is a metric $g^F$ on $F$ such that 
the scalar curvature of $g^F$ is positive. Under these two assumptions, 
A. Connes [Co] proved that one still has $\langle \widehat{A}(TM), [M]\rangle=0$. 
Note that here $M$ needs not  be spin and thus 
$\langle \widehat{A}(TM), [M]\rangle$ may not be an integer. Clearly this extends 
the Lichnerowicz 
theorem which corresponds to the special case $F=TM$.

The proof given by Connes in [Co] is highly 
noncommutative. It uses, besides an 
application of the Lichnerowicz formula along the leaves, 
 the longitudinal index theorem of
 Connes-Skandalis for foliated manifolds [CoS] as well as the cyclic 
cohomology techniques, see [Co].

The present paper arose from an attempt to find a direct geometric
 proof of Connes' vanishing theorem. As we will see, at least for the 
almost Riemannian foliations studied in [Co], Connes'
vanishing theorem can be directly deduced from a 
Lichnerowicz type formula and the Atiyah-Singer index theorem.

To make the use of the Lichnerowicz type formula to the 
foliated manifolds possible, 
we use the adiabatic limit procedure to blow up the metric along 
the transversal direction. Some geometric behaviors of the foliation 
under adiabatic limits are also examined. In particular, we find that 
the Bott connection [Bo] appears naturally in the limit. Then we construct
a class of natural elliptic operators which we called the sub-Dirac operators
associated to the foliation. The indices of such operators contain 
many important geometric information about the foliation. The 
method we use has more potentials. Actually we can prove certain 
new vanishing theorems for the foliated manifolds.

This paper is organized as follows. In \S 1 we discuss the geometry 
of a foliated manifold under the adiabatic limits. In \S 2, by 
constructing the sub-Dirac operators, we give a
direct geometric proof of 
the Connes vanishing theorem for what we call almost Riemannian foliations.
In \S 3, we briefly describe some new vanshing theorems that can be 
proved by using our method. The final section is an appendix in which we make some remarks
on the case of general foliations. 

$\ $

{\small {\bf Acknowledgements.} Part of this work was done while the first 
author was visiting the Nankai Institute of Mathematics in Tianjin 
during the summer of 1995. He would like to thank 
Prof. S. T. Yau who first proposed to use the idea of adiabatic limits 
to deal with Connes' theorem many years ago. Both of the authors wish 
to thank the Chinese National Science Foundation for its financial 
support. The first author was also partially supported by the NSF and the Sloan
Foundation, 
and the second author by the Ministry of Education of China and the
Qiu Shi Foundation. }

$\ $

{\bf \S 1. Adiabatic limits and foliations}

$\ $

The purpose of this section is to establish some formulas concering the adiabatic
limit of foliations.

Let $(M, F)$ be a  closed
foliation, that is, $F$ is an 
integrable subbundle of $TM$. Let $g^F$ be a metric on $F$. Let 
$g^{TM}$ be a metric on $TM$ which restricted to $g^F$ on $F$. 
Let $F^\bot$ be the orthogonal complement of $F$ in $TM$ 
with respect to $g^{TM}$. Then we have the following orthognal splittings,
$$TM=F\oplus F^\bot,$$
$$g^{TM}=g^F\oplus g^{F^\bot},\eqno(1.1)$$
where $g^{F^\bot}$ is the restriction of $g^{TM}$ to $F^\bot$.

It is clear that we can and we will make the identification that 
$$TM/F=F^\bot.\eqno(1.2)$$
Let $p,\ p^\bot$ be the orthogonal projection from $TM$ to $F$,
 $F^\bot$ respectively. Let $\nabla^{TM}$ be the Levi-Civita 
connection of $g^{TM}$ and $\nabla^F$ (resp. $\nabla^{F^\bot}$) be 
the restriction of $\nabla^{TM}$ to $F$ (resp. $F^{\bot}$). That is,    
$$\nabla^F= p\nabla^{TM}p,$$
$$\nabla^{F^\bot}= p^\bot \nabla^{TM}p^\bot.\eqno(1.3)$$

Now for any $\varepsilon>0$, let $g_\varepsilon^{TM}$ be the metric 
$$g^{TM,\varepsilon}=g^F\oplus {1\over \varepsilon} g^{F^\bot} .\eqno(1.4)$$ 
Let $\nabla^{TM,\varepsilon }$ be the Levi-Civita connection of 
$g^{TM,\varepsilon}$. Let $\nabla^{F,\varepsilon }$ (resp. 
$\nabla^{F^\bot,\varepsilon }$) be the restriction of 
$\nabla^{TM,\varepsilon }$ to $F$ (resp. $F^\bot$). 
We will examine the behavior of $\nabla^{TM,\varepsilon }$ as 
$\varepsilon\rightarrow 0$. 

The process of taking the limit 
$\varepsilon\rightarrow 0$ is called adiabatic limit. 

The standard formula for Levi-Civita connection gives us the 
following sets of formulas,
$$  \nabla^{F,\varepsilon }=\nabla^F,$$  
$$p\nabla_X^{TM,\varepsilon}p^\bot=p\nabla_X^{TM}p^\bot,\ 
\mbox{for} \ X\in \Gamma(F),\eqno(1.5)$$
and 
$$\left\langle\nabla_V^{TM,\varepsilon}U, X\right\rangle
=\left\langle\nabla_V^{TM}U, X\right\rangle
-{1\over 2}\left\langle X,\nabla_V^{TM}U+ \nabla_U^{TM}V\right\rangle
+ {1\over 2\varepsilon }\left\langle X, \nabla^{TM}_V U+\nabla^{TM}_UV\right\rangle, 
\eqno(1.6)$$
for $X\in \Gamma (F), \ U,\ V \in \Gamma(F^\bot)$. 
Furthermore,
$$p^\bot\nabla_X^{TM,\varepsilon}p=\varepsilon p^\bot\nabla_X^{TM}p,\ 
\mbox{for}\  X\in \Gamma(F),$$
$$\left\langle \nabla_V^{TM,\varepsilon}Y, U\right\rangle=-{1\over 2}\left\langle
Y,\nabla^{TM}_V U+\nabla^{TM}_UV\right\rangle
+{\varepsilon\over 2}\langle Y,[U,V]\rangle\eqno(1.7)$$
and 
$$\nabla_V^{F^\bot,\varepsilon}=\nabla_V^{F^\bot}, $$
$$\left\langle \nabla_X^{F^\bot,\varepsilon}U, V\right\rangle=
\left\langle [X,U], V\right\rangle -
{1\over 2}\left\langle X,\nabla^{TM}_V U+\nabla^{TM}_UV\right\rangle
-{\varepsilon\over 2}\langle X,[U,V]\rangle \eqno (1.8)$$
$\mbox{for} \ X\in \Gamma(F), \ U,\ V \in \Gamma(F^\bot).$ 

Let $L$ be a leave of $F$, then $F^\bot$ is a flat bundle 
along $L$, carrying with the canonical Bott connection [Bo]
$$\dot{\nabla}_X^L=p^\bot[X,U], \ X\in \Gamma(TL),\ U\in 
\Gamma(F^\bot).\eqno(1.9)$$

Following [BZ], set 
$$\omega^L=\left(g^{F^\bot}\right)^{-1}\dot{\nabla}^L g^{F^\bot}=
\dot{\nabla}^{L*}-\dot{\nabla}^L ,\eqno(1.10)$$ 
where $\dot{\nabla}^{L*}$ is the dual of $\dot{\nabla}^{L}$ 
with respect to $g^{F^\bot}$. Then 
$$\widetilde{\nabla}^L=\dot{\nabla}^L+{1\over 2}\omega^L \eqno(1.11)$$ 
is the natural unitary connection on $F^\bot|_L$ associated to $\dot{\nabla}^L$.

We can now state the main result of this section as follows.

$\ $ 

{\bf Theorem 1.1.} {\em  Along each leave $L$, the following identity holds,}
$$\lim_{\varepsilon\rightarrow 0}\left.\nabla^{F^\bot,\varepsilon}\right|_L = 
\widetilde{\nabla}^L.\eqno(1.12)$$ 

{\em Proof.}  For any $X\in \Gamma(TL)$, $U,\ V\in \Gamma(F^\bot)$, one has 
$$\omega^L(X)(U,V)=\left\langle \dot{\nabla}^{L*}_X U, V\right\rangle-
\left\langle\dot{\nabla}^L_X U, 
V\right\rangle$$
$$=-\left\langle U, \dot{\nabla}^L_XV\right\rangle
-\left\langle \dot{\nabla}^L_X U,V\right\rangle+X\langle U, V\rangle$$
$$=-\langle U,[X,V]\rangle -\langle [X, U], V\rangle+X\langle U,V\rangle$$
$$=-\left\langle U, \nabla^{TM}_XV-\nabla^{TM}_VX\right\rangle 
-\left\langle \nabla^{TM}_XU-\nabla^{TM}_UX, V\right\rangle +X\langle U,V\rangle $$
$$=-\left\langle\nabla^{TM}_UV,X\right\rangle-
\left\langle\nabla^{TM}_VU,X\right\rangle
-\left\langle U,\nabla^{TM}_XV\right\rangle-
\left\langle\nabla^{TM}_XU,V\right\rangle+X\langle U,V\rangle.
\eqno(1.13)$$

Note that the last three terms cancel. So (1.12) follows directly 
from (1.8), (1.9), (1.11) and (1.13). $\Box$

$\ $

{\bf Remark 1.2.} Conversely, one sees from Theorem 1.1 that the Bott 
connection shows up naturally from the above adiabatic limit procedure.

$\ $

{\bf \S 2. The Connes vanishing theorem for almost Riemannian foliations}

$\ $

In this section we present a direct geometric proof
of the Connes vanishing theorem for a special class 
of foliated manifolds which we call almost Riemannian foliations. 

This section is organized as follows. In a), we define what we call  
almost Riemannian foliations. In b) we construct a class of elliptic operators 
specially defined for our purpose and prove the corresponding Lichnerowicz 
formula. In c) we prove the Connes vanishing theorem for almost 
Riemannian foliations.

$\ $

{\bf a). Almost Riemannian foliations}

$\ $

We use the same assumptions and notations as in \S 1. First note that the 
$1$-form $\omega^L$ defined in
 (1.10) does not depend on the metric $g^F$. Also recall that, if $\omega=0$, 
then $(M,F,g^{TM})$ is a Riemannian foliation. This 
motivates the following definition.

$\ $

{\bf Definition 2.1.} Let $g^F$ be a metric on $F$. If there is a series of 
metrics $g^{TM}_\varepsilon $, $\varepsilon>0$, such that as 
$\varepsilon\rightarrow 0$, the corresponding one form defined in (1.10)
verifies that
$$\omega_\varepsilon\rightarrow 0,\eqno(2.1)$$ 
then we say that $(M, F, g^F)$ admits an almost Riemannian structure.

$\ $

{\bf Remark 2.2.} The equation (2.1) should be interpreted as follows: as 
$\varepsilon \rightarrow 0$, for any $X\in \Gamma(F)$, one has 
$$ |\omega_\varepsilon(X)|_{g^{TM}_\varepsilon}\rightarrow 0.\eqno
(2.2)$$

$\ $

{\bf b). A class of elliptic operators associated to Spin subbundles of 
the tangent bundle}

$\ $

{}From now on we make the special assumption that $F$ is oriented,
spin and carries 
a fixed spin structure. We also assume that $F^\bot$ is oriented and that 
both $p=\dim F$ and $q=\dim F^\bot$ are even.

Let $S(F)$ be the bundle of spinors associated to $(F, g^F)$. For any 
$X\in \Gamma(F)$, denote by $c(X)$ the Clifford action of $X$ on $S(F)$. 
Since $p=\dim F$ is even, we have the splitting 
$$S(F)=S_+(F)\oplus S_-(F)\eqno(2.3)$$
and $c(X)$ exchanges $S_\pm(F)$.

Let $\Lambda(F^{\bot,*})$ be the exterior algebra bundle of $F^\bot$. 
Then $\Lambda(F^{\bot,*})$ carries a canonically induced metric 
$g^{\Lambda(F^{\bot,*})}$ from $g^{F^\bot}$.
For any $U\in \Gamma(F^\bot)$, 
let $U^*\in \Gamma(F^{\bot, *})$ be the corresponding dual of 
$U$ with respect to $g^{F^\bot}$.

Now for $U\in \Gamma(F^\bot)$, set 
$$ c(U)=U^*\land -i_U, \ \widehat{c}(U)=U^*\land+i_U,\eqno(2.4)$$ 
where $U^*\land $ and $i_U$ are the exterior and inner multiplications by 
$U^*$ and $U$ on $\Lambda(F^{\bot,*})$ respecyively. One has the following 
obvious identities,
$$c(U)c(V)+c(V)c(U)=-2\langle U,V\rangle_{g^{F^\bot}},$$ 
 $$\widehat{c}(U)\widehat{c}(V)+\widehat{c}(V)\widehat{c}(U)
=2\langle U,V\rangle_{g^{F^\bot}},$$
$$c(U)\widehat{c}(V)+\widehat{c}(V)c(U)=0  \eqno(2.5)$$
for $  U,\ V\in \Gamma (F^\bot)$.

Let $h_1,\cdots, h_q$ be an oriented local orthonormal basis of $F^\bot$. 
Set 
$$\tau\left(F^\bot, g^{F^\bot}\right)=\left(
{1\over \sqrt{-1}}\right)^{q(q+1)\over 2}c(h_1)\cdots c(h_q).
\eqno(2.6)$$
Then  
$$\tau\left(F^\bot, g^{F^\bot}\right)^2={\rm Id}_{\Lambda(F^{\bot,*})}.
\eqno(2.7)$$

Denote 
$$\Lambda_\pm\left(F^{\bot,*}\right)=\left\{ h\in \Lambda\left(F^{\bot,*}\right):\ 
\tau\left(F^\bot, 
g^{F^\bot}\right)h=\pm h \right\}.\eqno(2.8)$$
Then $\Lambda_\pm(F^{\bot,*})$ are sub-bundles of $\Lambda(F^{\bot,*})$. 
Also, since $q$ is even, one verifies that for any $h\in \Gamma(F^\bot)$, 
$c(h)$ anticommutes with $\tau$, while $\widehat{c}(h)$ commutes with $\tau$. 
Thus $c(h)$ exchanges $\Lambda_\pm(F^{\bot,*})$.

We will view both vector bundles 
$$S(F)=S_+(F)\oplus S_-(F)\eqno(2.9)$$ 
and 
$$\Lambda\left(F^{\bot,*}\right)=\Lambda_+\left(F^{\bot,*}\right)\oplus \Lambda_-
\left(F^{\bot,*}\right)\eqno(2.10)$$ 
as super vector bundles. 
Their ${\bf Z}_2$ graded tensor product is given by 
$$S(F)\widehat{\otimes}\Lambda(F^{\bot,*})=\left[S_+(F)\otimes\Lambda_+
\left(F^{\bot,*}\right)
\oplus S_-(F)\otimes\Lambda_-\left(F^{\bot,*}\right)\right]$$
$$\bigoplus\left[S_+(F)\otimes\Lambda_-\left(F^{\bot,*}\right)
\oplus S_-(F)\otimes\Lambda_+\left(F^{\bot,*}\right)\right].\eqno(2.11)$$
For $X\in \Gamma(F)$, $U\in \Gamma(F^\bot)$, the operators $c(X)$, 
$c(U)$ and $\widehat{c}(U)$ extends naturally to 
$S(F)\widehat{\otimes}\Lambda(F^{\bot,*})$.

The connections $\nabla^F$, $\nabla^{F^\bot}$ lift to $S(F)$ and 
$ \Lambda(F^{\bot,*})$ naturally , and preserve the splittings 
(2.9) and (2.10). We write them as
$$\nabla^{S(F)}=\nabla^{S_+(F)}\oplus\nabla^{S_-(F)},$$
$$\nabla^{\Lambda(F^{\bot,*})}=\nabla^{\Lambda_+(F^{\bot,*})}
\oplus\nabla^{\Lambda_-(F^{\bot,*})}.\eqno(2.12)$$
Then $S(F)\widehat{\otimes}\Lambda(F^{\bot,*})$ carries the induced 
tensor product connection
$$\nabla^{S(F)\widehat{\otimes}\Lambda(F^{\bot,*})}=\nabla^{S(F)}\otimes 
{\rm Id}_{\Lambda(F^{\bot,*})}+{\rm Id}_{S(F)}\otimes\nabla^{\Lambda(F^{\bot,*})} .
\eqno(2.13)$$
And similarly for $S_\pm(F)\widehat{\otimes}\Lambda_\pm(F^{\bot,*})$.

Let $S\in \Omega^1(T^*M)\otimes \Gamma(\mbox{End}(TM))$ be defined by 
$$\nabla^{TM}=\nabla^{F}+\nabla^{F^\bot}+S.\eqno(2.14)$$ 
Then for any $X\in \Gamma(TM)$, 
$S(X)$ exchanges $\Gamma(F)$ and $\Gamma(F^\bot)$ and is skew-adjoint with
 respect to $g^{TM}$.

For any vector bundle $E$ over $M$. By an integral polynomial of $E$ we 
will mean a vector bundle $\phi(E)$ which is a polynomial in the exterior and
 symmetric powers of $E$ with integral coefficients.

Let $\phi(F^\bot)$ be an integral polynomial of $F^\bot$, then $\phi(F^\bot)$ 
carries a naturally induced metric $g^{\phi(F^\bot)}$ from $g^{F^\bot}$ and also 
a naturally induced Hermitian connection $\nabla^{\phi(F^\bot)}$ induced from 
$\nabla^{F^\bot}$.

Our main concern will be on the $Z_2$-graded vector bundle 
$$\left(S(F)\widehat{\otimes}\Lambda\left(F^{\bot,*}\right)
\right)\otimes\phi(F^\bot)\eqno(2.15)$$ 
which is 
$$\left[S_+(F)\otimes\Lambda_+\left(F^{\bot,*}\right)\otimes\phi\left(F^\bot\right) 
\oplus S_-(F)\otimes\Lambda_-\left(F^{\bot,*}\right)\otimes\phi
\left(F^\bot\right)\right]$$
$$\bigoplus\left[S_+(F)\otimes\Lambda_-\left(F^{\bot,*}\right)\otimes
\phi\left(F^\bot\right) 
\oplus S_+(F)\otimes\Lambda_-\left(F^{\bot,*}\right)
\otimes\phi\left(F^\bot\right)\right].$$

The Clifford actions $c(X)$, $c(U)$ and $\widehat{c}(U)$ for $X\in\Gamma(F)$, 
$U\in \Gamma(F^\bot)$ extends further to these bundles by acting as identity
 on $\phi(F^\bot)$.

We can also form the tensor product metric on the new bundles as well as 
the tensor product connection
$$\nabla^{(S(F)\widehat{\otimes}\Lambda(F^{\bot,*}))\otimes \phi(F^\bot)}
=\nabla^{S(F)\widehat{\otimes}\Lambda(F^{\bot,*})}\otimes {\rm Id}_{\phi(F^\bot)}+
{\rm Id}_{S(F)\widehat{\otimes}\Lambda(F^{\bot,*})}
\otimes \nabla^{\phi(F^\bot)},\eqno(2.16)$$
and similarly for the $\pm$ subbundles.

Now let $\{ f_i\}^p_{i=1}$ be an oriented orthonormal 
basis of $F$. Recall 
that $\{ h_s\}^q_{s=1}$ is an orthonormal basis of $F^\bot$. The elliptic 
operators which are the main concern of this subsection can be defined as 
follows. It is introduced mainly for the reason that the vector bundle
$F^\bot$ might well be non-spin.

$\ $

{\bf Definittion 2.2.} Let $D_{F, \phi(F^\bot)}$ be the 
operator mapping from 
$\Gamma(S(F)\widehat{\otimes}\Lambda(F^{\bot,*})\otimes \phi(F^\bot))$ 
to itself defined by  
$$D_{F, \phi(F^\bot)}= \sum^p_{i=1}c(f_i)\nabla_{f_i}^{(S(F)
\widehat{\otimes}\Lambda(F^{\bot,*}))\otimes \phi(F^\bot)}
 +\sum^q_{s=1}c(h_s)
\nabla_{h_s }^{(S(F)\widehat{\otimes}\Lambda(F^{\bot,*}))\otimes \phi(F^\bot)}$$
$$+{1\over 2}\sum_{i,j=1}^p\sum^q_{s=1}\langle
 S(f_i)f_j, h_s\rangle c(f_i)c(f_j)c(h_s)
 +{1\over 2}\sum_{s,t=1}^q\sum^p_{i=1}
\langle S(h_s)h_t,f_i\rangle c(h_s)c(h_t)c(f_i).\eqno(2.17)$$

It is easy to verify that $D_{F, \phi(F^\bot)}$ is a first order 
formally self adjoint elliptic differential operator. Furthermore it 
anticommutes with the ${\bf Z}_2$ grading operator of the super vector bundle 
$(S(F)\widehat{\otimes}\Lambda(F^{\bot,*}))\otimes \phi(F^\bot)$.

Let $D_{+,F, \phi(F^\bot)}$ (resp. $D_{-,F, \phi(F^\bot)}$) be the 
restriction of $D_{F, \phi(F^\bot)}$ to the even (resp. odd) subbundle 
of $(S(F)\widehat{\otimes}\Lambda(F^{\bot,*}))\otimes \phi(F^\bot)$. 
Then one has 
$$D_{+,F, \phi(F^\bot)}^*=D_{-,F, \phi(F^\bot)}.\eqno(2.18)$$

Let $\Delta^{F,\phi(F^\bot)}$ be the Bochner Laplacian defined by 
$$\Delta^{F,\phi(F^\bot)}= 
 \sum^p_{i=1}\left(\nabla_{f_i}^{(S(F)\widehat{\otimes}\Lambda(F^{\bot,*}))
\otimes \phi(F^\bot)}+{1\over 2}
\sum_{j=1}^p\sum_{s=1}^q\left\langle S(f_i)f_j, h_s\right\rangle
c(f_j)c(h_s)\right)^2$$
$$+\sum^q_{s=1}\left(\nabla_{h_s}^{(S(F)\widehat{\otimes}
\Lambda(F^{\bot,*}))\otimes \phi(F^\bot)}+{1\over 2}\sum_{t=1}^q
\sum_{j=1}^p\langle S(h_s)h_t,f_j\rangle c(h_t)c(f_j)\right)^2$$
$$-\left(\nabla^{(S(F)\widehat{\otimes}\Lambda(F^{\bot,*}))\otimes 
\phi(F^\bot)}_{\Sigma_{i=1}^p\nabla^{TM}_{f_i}f_i}+{1\over 2}\sum_{j=1}^p
\sum_{s=1}^q \left\langle S\left(\sum_{i=1}^p\nabla^{TM}_{f_i}f_i\right)f_j, 
h_s\right\rangle c(f_j)c(h_s)\right)$$
$$-\left(\nabla^{(S(F)\widehat{\otimes}\Lambda(F^{\bot,*}))\otimes 
\phi(F^\bot)}_{\Sigma_{s=1}^q\nabla^{TM}_{h_s}h_s}+{1\over 2}\sum_{t=1}^q
\sum_{j=1}^p\left\langle S\left(\sum_{s=1}^q \nabla^{TM}_{h_s}h_s\right)h_t, 
f_j\right\rangle
c(h_t)c(f_j)\right).\eqno(2.19)$$

Let $k_{TM}$ be the scalar curvature of the metric $g^{TM}$, 
let $R^{F^\bot}$ be the curvature tensor of 
$\nabla^{F^\bot}$. Let $R^{\phi(F^\bot)}$ be the curvature of 
$\nabla^{\phi(F^\bot)}$. Then we can state the corresponding 
Lichnerowicz type formula for  $D_{F, \phi(F^\bot)}$ as follows.

$\ $

{\bf Theorem 2.3.} {\em The following identity holds,
$${D^2_{F,\phi(F^\bot)}}= -\Delta^{F,\phi(F^\bot)}
+ {1\over 2}\sum_{i,j=1}^pc(f_i)c(f_j)R^{\phi(F^\bot)}(f_i,f_j)$$
$$+\sum_{i=1}^p\sum_{s=1}^qc(f_i)c(h_s)R^{\phi(F^\bot)}
(f_i,h_s)
+{1\over 2}\sum_{s,t=1}^qc(h_s)c(h_t)R^{\phi(F^\bot)}(h_s,h_t)$$
$$+{k_{TM}\over 4}
+{1\over 4}\sum_{i=1}^p\sum_{r,s,t=1}^q\left\langle R^{F^\bot}
(f_i,h_r)h_t,h_s\right\rangle
c(f_i)c(h_r)\widehat{c}(h_s)\widehat{c}(h_t)$$
$$+{1\over 8}\sum_{i,j=1}^p\sum_{s,t=1}^q\left\langle R^{F^\bot}
(f_i,f_j)h_t,h_s\right\rangle
c(f_i)c(f_j)\widehat{c}(h_s)\widehat{c}(h_t)$$
$$+{1\over 8}\sum_{s,t,r,l=1}^q\left\langle R^{F^\bot}(h_r,h_l)h_t,h_s
\right\rangle c(h_r)c(h_l)
\widehat{c}(h_s)\widehat{c}(h_t).\eqno(2.20) $$}

{\it Proof.} We first assume that $F^\bot$, and thus $TM$ also, is spin. 
Then $D_{F, \phi(F^\bot)}$ is just the standard Dirac operator on 
$M$ twisted by $S(F^\bot)\otimes \phi(F^\bot)$ which carries the 
canonical connections induced from $\nabla^{F^\bot}$,
$$D_{F, \phi(F^\bot)}:  \Gamma\left(S(TM)\otimes S\left(F^\bot\right)\otimes 
\phi\left(F^\bot\right)\right)
\longrightarrow \Gamma\left(S(TM)\otimes 
S\left(F^\bot\right)\otimes 
\phi\left(F^\bot\right)\right).\eqno(2.21)$$
With this observation, (2,20) is a direct corollary of the usual 
Lichnerowicz formula [L] for twisted Dirac operators.

Now observe that (2.20) is a local formula, and that
locally we can always assume $TM$ is spin. 
Thus the formula is proved. The interested reader may also proceed to 
verify (2.20) by a direct computation.
$\Box$

$\ $

{\bf Remark 2.4.}  The definition of $D_{F,\phi(F^\bot)}$ does not use 
the condition that $F$ is integrable.

$\ $

{\bf Remark 2.5.}  While locally $D_{F,\phi(F^\bot)}$ can be seen as a 
twisted Dirac operator, the key point here is that its definition 
relies {\em only} on the spin structure of $F$.

$\ $

{\bf c). Almost Riemannian foliations and the Connes vanishing theorem}

$\ $

Recall that $(M,F)$ is a compact foliated manifold and we assume 
that $F$ is oriented, spin and both $p=\dim F$ and $q=\dim F^\bot$ 
are even, and that $TM/F$ is oriented.

For any $x\in M$, let $L_x$ be the leave passing through $x$. Then 
$g^{TF}$ restricts to a Riemannian metric on $L_x$. So it determines a 
scalar curvuture $k_F(x)$.

We now make the basic assumption that $k_F$ is positive on $M$. Since
 $M$ is compact, there is then a positive number $\delta>0$, such that 
$$k_F(x)>\delta\ \ \ \mbox{for any }\ x\in M.\eqno(2.22)$$

The second important assumption for this subsection is that $(M, F, g^F)$ 
is an almost Riemannian foliated manifold in the sense of Definition 2.1.

We fix temporarily a splitting 
$$TM=F\oplus F^\bot.\eqno(2.23)$$ 
Given a metric $g^{F^\bot}$ on $F^\bot$, recall that we have  defined  
$$\omega=\left(g^{F^\bot}\right)^{-1}\dot{\nabla}g^{F^\bot}\eqno(2.24)$$ 
along the leaves of $F$, or we can say, along $F$.

By [BZ], one sees  that
$$\dot{\nabla}^2=-{1\over 4} \omega^2\eqno (2.25)$$
along $F$.

Let $\phi(F^\bot)$ be an integral polynomial of $F^\bot$. Let 
$\dot{\nabla}_\phi$ be the lift of $\dot{\nabla}$ to $\phi(F^\bot)$ 
along $F$. Then $\dot{\nabla}_\phi$ is also a flat connection on
 $\phi(F^\bot)$ along $F$.

Let $\omega_\phi$ be the lift of $\omega$ to $ \phi(F^\bot)$. It 
is easy to see that $\omega_\phi$ can be represented as linear
 combinations of the powers of $\omega$. Inparticular, if  
$$||\omega||:=\sup_{x\in M} \{ |\omega(X)|:X\in \Gamma(F), |X|\leq 1\} <1,\eqno(2.26)$$ 
where $|\cdot|$ is the norm 
with respect to the metric $g^{TM}$, then there exists a constant 
$C_\phi>0$ such that for any $X\in \Gamma(F)$ with  $|X|\leq 1$,
$$|\omega_\phi(X)|\leq C_\phi ||\omega||.\eqno(2.27)$$

Also, similar to (2.25), we have 
$$\left(\dot{\nabla}_\phi\right)^2=-{1\over 4} \omega^2_\phi\eqno
(2.28)$$ 
along $F$.

Now for any $\varepsilon>0$, let $g^{TM,\varepsilon}$ be the metric 
$$g^{TM,\varepsilon}=g^F\oplus {1\over \varepsilon} g^{F^\bot}.\eqno
(2.29)$$
Let $D_{F, \phi(F^\bot), \varepsilon}$ be the elliptic operator as 
constructed in last subsection but for the metric $g^{TM,\varepsilon}$.
 Let $\Delta^{F, \phi(F^\bot), \varepsilon}$ be the corresponding Bochner 
Laplacian. We will use suitable subscript or superscript on the 
corrresponding gemetric quantities to indicate that they are with 
respect to $g^{TM,\varepsilon}$. 

We now examine the behaviour of
 $D_{F, \phi(F^\bot), \varepsilon}^2$ when $\varepsilon\rightarrow 0$.

Recall that $\{f_i\}^p_{i=1}$, $\{ h_s\}^q_{s=1}$ constitute an
 orthonormal basis of $g^{TM}=g^F\oplus g^{F^\bot}$. Thus 
$\{f_i\}^p_{i=1}$, $\{ \sqrt{ \varepsilon }h_s\}^q_{s=1}$ is an 
orthonormal basis for $g^{TM,\varepsilon}$.

By applying Theorem 2.3 to this situation, we get 
$${D^2_{F,\phi(F^\bot),\varepsilon }}= -\Delta^{F,\phi(F^\bot),\varepsilon }
+ {1\over 2}\sum_{i,j=1}^pc(f_i)c(f_j)R^{\phi(F^\bot),\varepsilon }(f_i,f_j)$$
$$+\sqrt{\varepsilon }\sum_{i=1}^p\sum_{s=1}^qc(f_i)c(h_s)
R^{\phi(F^\bot),\varepsilon }
(f_i,h_s)
+{\varepsilon \over 2}\sum_{s,t=1}^qc(h_s)c(h_t)R^{\phi(F^\bot),\varepsilon }(h_s,h_t)$$
$$+{k_{TM,\varepsilon }\over 4}
+{\sqrt{\varepsilon }\over 4}\sum_{i=1}^p\sum_{r,s,t=1}^q\left\langle 
R^{F^\bot,\varepsilon }
(f_i,h_r)h_t,h_s\right\rangle
c(f_i)c(h_r)\widehat{c}(h_s)\widehat{c}(h_t) $$
$$+{1\over 8}\sum_{i,j=1}^p\sum_{s,t=1}^q\left\langle R^{F^\bot,\varepsilon }
(f_i,f_j)h_t,h_s\right\rangle
c(f_i)c(f_j)\widehat{c}(h_s)\widehat{c}(h_t)$$
$$+{\varepsilon \over 8}\sum_{s,t,r,l=1}^q\left\langle 
R^{F^\bot,\varepsilon }(h_r,h_l)h_t,h_s
\right\rangle c(h_r)c(h_l)
\widehat{c}(h_s)\widehat{c}(h_t).\eqno(2.30) $$

We first examine the behaviour of $k_{TM,\varepsilon}$ as 
$\varepsilon\rightarrow 0$. 

By definition, one has 
$$k_{TM,\varepsilon}=\sum_{i,j=1}^p\left\langle R^{TM,\varepsilon}
(f_i, f_j)f_i, f_j\right\rangle
 +\sum_{s,t=1}^q\varepsilon\left\langle R^{TM,\varepsilon}(h_s, h_t)h_s, h_t
\right\rangle 
    +\sum_{i=1}^p\sum_{s=1}^q\left\langle R^{TM,\varepsilon}(f_i, h_s)f_i,h_s
\right\rangle,\eqno(2.31) $$
where $R^{TM,\varepsilon}$ is the curvature of $\nabla^{TM,\varepsilon}$, 
the Levi-Civita connection of $g^{TM,\varepsilon}$.

Recall that $p, \ p^\bot$ are the orthogornol projections, with respect 
to $g^{TM}$, from $TM$ to $F, F^\bot$ respectively. From (1.5) to (1.8), 
we deduce that 
$$\left\langle R^{TM,\varepsilon}(f_i, f_j)f_i, f_j\right\rangle=
\left\langle R^{F}(f_i, f_j)f_i, f_j\right\rangle
+\varepsilon\left\langle \nabla^{TM}_{f_j}f_j, 
p^\bot\nabla^{TM}_{f_i}f_i\right\rangle
-\varepsilon\left\langle\nabla^{TM}_{f_i}f_j, 
p^\bot\nabla^{TM}_{f_j}f_i\right\rangle,
\eqno(2.32) $$ 
and, for $X\in \Gamma(F), \ U\in \gamma(F^\bot)$ that
$$\left\langle R^{TM,\varepsilon}(X, U)X, U\right\rangle=
\varepsilon\left\langle\nabla^{TM}_Xp\nabla^{TM}_UX, U\right\rangle
-\left\langle\nabla^{TM}_U\left(p\nabla^{TM,\varepsilon}_XX\right),U\right\rangle$$
$$-\varepsilon\left\langle\nabla^{TM}_{p[X,U]}X, U\right\rangle
-\left\langle\nabla^{TM}_{p^\bot[X,U]}X, U\right\rangle+{1\over 2}
\left\langle X,\left[U,p^\bot[X,U]\right]\right\rangle$$
$$ -{\varepsilon\over 2}\left\langle X,\left[U,p^\bot[X,U]\right]\right\rangle  
-{1\over 2}\left\langle X,\nabla^{TM}_{p^\bot\nabla^{TM,\varepsilon}_UX}U+
\nabla^{TM}_U\left(p^\bot\nabla^{TM,\varepsilon}_UX\right)\right\rangle$$
$$+\left\langle \left[X, p^\bot\nabla^{TM,\varepsilon}_UX\right], U\right\rangle
+{\varepsilon\over 2}
 \left\langle X,\left[U,p^\bot\nabla^{TM,\varepsilon}_UX\right]
\right\rangle-
\left\langle \nabla^{TM}_U\left(p^\bot\nabla^{TM,\varepsilon}_XX
\right), U\right\rangle$$
$$=-\varepsilon \left\langle p\nabla^{TM}_UX,\nabla^{TM}_XU\right\rangle -
\varepsilon \left\langle \nabla^{TM}_{p[X,U]}X, U\right\rangle
-{\varepsilon\over 2}\left\langle X,\left[U,p^\bot[X,U]\right]\right\rangle$$
$$+{\varepsilon\over 2}
\left\langle X,p^\bot\left(\nabla^{TM,\varepsilon}_UX\right)\right\rangle
-\varepsilon\left\langle\nabla^{TM}_U\left(p^\bot
\nabla^{TM}_XX\right),
U\right\rangle+\left\langle p\nabla^{TM}_XX,\nabla^{TM}_UU\right\rangle$$
$$+{1\over 2}\left\langle X,\nabla^{TM}_U\left(p^\bot[X,U]\right)+
\nabla^{TM}_{p^\bot[X,U]}U\right\rangle +
\left\langle \left[X,p^\bot\nabla^{TM,\varepsilon}_UX\right], U\right\rangle$$ 
$$-{1\over 2}\left\langle X,\nabla^{TM}_{p^\bot
\nabla^{TM,\varepsilon}_UX}U
+\nabla^{TM}_U\left(p^\bot\nabla^{TM,\varepsilon}_UX\right)\right\rangle,
\eqno(2.33)$$
while for $U, \, V\in \Gamma(F^\bot)$, we have 
$$\left\langle R^{TM,\varepsilon}(U, V) U,V\right\rangle
=-{1\over 2}\left\langle p\nabla^{TM,\varepsilon}_VU,
\nabla^{TM}_UV+\nabla^{TM}_VU\right\rangle 
+\left\langle\nabla^{TM}_{p[U,V]}U,V\right\rangle$$
$$ +{\varepsilon\over 2}
\left\langle p\nabla^{TM,\varepsilon}_VU,[V,U]\right\rangle 
+{1\over 2}\left\langle p\nabla^{TM,\varepsilon}_VU,\nabla^{TM}_UV+
\nabla^{TM}_VU\right\rangle$$
$$ - {\varepsilon\over 2}
\left\langle p\nabla^{TM,\varepsilon}_UV,[V,U]\right\rangle 
+\left\langle R^{F^\bot}(U,V)U,V\right\rangle-
\left\langle\nabla^{TM,\varepsilon}_{p[U,V]}U,V\right\rangle$$
$$=\left\langle R^{F^\bot}(U,V)U,V\right\rangle+{1\over 2}
\left\langle [U,V],p\left(\nabla^{TM}_UV+\nabla^{TM}_VU\right)\right\rangle
-\left({1\over 2}-\varepsilon\right)\left\langle 
p[U,V],[U,V]\right\rangle.\eqno(2.34)$$

On the other hand, by (1.7), one has 
$$p^\bot\nabla^{TM,\varepsilon}_UX=\sum_{s=1}^q
\left\langle\nabla^{TM,\varepsilon}_UX,h_s\right\rangle
h_s $$
$$=-\sum_{s=1}^q \left\{ {1\over 2} 
\left\langle X,\nabla^{TM}_Uh_s+\nabla^{TM}_{h_s}U\right\rangle
h_s+{\varepsilon\over 2}\left\langle X, [U,h_s]\right\rangle h_s\right\}$$
$$={1\over 2} \sum_{s=1}^q
\omega(X)(U,h_s)h_s
-{\varepsilon\over 2}\sum_{s=1}^q\left\langle X,[U,h_s]\right\rangle h_s.
\eqno(2.35)$$
Here we have used the notation of \S 1 and also formula (1.13). 
Also since the Bott connection (1.9) is flat, by (1.10) one
 verifies directly that (cf. [BZ])
$$\dot{\nabla}\omega=-{1\over 4} \omega^2, \eqno(2.36)$$ 
along $F$.

Equations  (2.35) and (2.36) can be used to control the term 
$$\left\langle \left[X, p^\bot\nabla^{TM,\varepsilon}_UX\right],
 U\right\rangle \eqno(2.37)$$ 
in (2.33). Combining with (1.13), (2.22) and $(2.32)$--$(2.36)$, 
one finds that there are positive constants $\varepsilon_0>0$, $C>0$ 
such that if $0<\varepsilon<\varepsilon_0$, then 
$$k_{TM,\varepsilon}>\delta-C||\omega||.\eqno(2.38)$$

Now examining the other curvarure terms appearing in (2.30). This part 
is easier, as we have the convergence formula (1.12).

Note that $\dot{\nabla}$ as well as $\dot{\nabla}_\phi$ are flat 
along $F$. By (1.11), (2.27) and (2.28), we know that whenever 
$||\omega||<1$, there exist positive constants 
$C', \ \varepsilon_0'$, such that if 
$0< \varepsilon< \varepsilon_0'$, then 

$$\left|{1\over 8}\sum_{i,j=1}^p\sum_{s,t=1}^q\left\langle R^{F^\bot, 
\varepsilon}(f_i,f_j)h_t,h_s)\right\rangle
c(f_i)c(f_j)\widehat{c}(h_s)\widehat{c}(h_t)\right.$$
$$+{\sqrt{\varepsilon }\over 4}\sum_{i=1}^p\sum_{r,s,t=1}^q\left\langle 
R^{F^\bot,\varepsilon }
(f_i,h_r)h_t,h_s\right\rangle
c(f_i)c(h_r)\widehat{c}(h_s)\widehat{c}(h_t) $$
$$+{\varepsilon\over 8}\sum_{s,t,r,l=1}^q
\left\langle R^{F^\bot, \varepsilon}(h_r,h_l)h_t,h_s\right\rangle
c(h_r)c(h_l)\widehat{c}(h_s)\widehat{c}(h_t)$$
$$+{1\over 2} \sum_{i,j=1}^pc(f_i)c(f_j)R^{\phi(F^\bot), 
\varepsilon}(f_i,f_j)+\sqrt{\varepsilon}\sum_{i=1}^p\sum_{s=1}^q
c(f_i)c(h_s)R^{\phi(F^\bot), \varepsilon}(f_i,h_s)$$
$$\left. +{\varepsilon\over 2}\sum_{s,t=1}^qc(h_s)c(h_t)
R^{\phi(F^\bot), \varepsilon}(h_s,h_t)\right|\leq C'||\omega||.
\eqno(2.39)$$

Now since we have assumed that $(M,F, g^F)$ is almost Riemannian, 
we can and we will choose $g^{TM}$ such that  
$$||\omega||<1,\ \ \ \left(C+C'\right)||\omega||<{\delta\over 8}.   \eqno(2.40)$$

{}From (2.30), (2.38), (2.39) and (2.40), we see that there is a 
positive constant $\varepsilon_0''$ such that when 
$0<\varepsilon<\varepsilon_0''$, we have 

$$D_{F,\phi(F^\bot),\varepsilon}^2>-\Delta^{F,\phi(F^\bot),
\varepsilon}+{\delta\over 8}>0,\eqno(2.41)$$
since $-\Delta^{F,\phi(F^\bot),\varepsilon}$ is clearly nonnegative.

Now we can prove the following result which is the Connes 
vanishing theorem for  almost Riemannian foliations.

$\ $

{\bf Theorem 2.6.} {\em Let $(M,F, g^F)$ be an almost Riemannian 
foliation, assume $M$ is compact and transversally oriented. 
If $F$ is spin and $k_F>0$, then 
$$\left\langle\widehat{A}(F)p\left(F^\bot\right),[M]\right\rangle =0,
\eqno(2.42)$$
where $p(F^\bot)$ is any Pontrjagin class of $F^\bot$.}

{\bf Proof:}  By (2.18), (2.41), one deduces that
$$\mbox{ind}\, D_{+,F,\phi(F^\bot),\varepsilon}=0 .\eqno(2.43)$$

Thus by applying the Atiyah-Singer index theorem [AS], or 
by a direct heat kernel evaluation, we get 
$$\mbox{ind}\, D_{+,F,\phi(F^\bot)}=\int_M\widehat{A}(F)L\left(F^\bot\right)\mbox{ch}
 \left(\phi\left(F^\bot\right)\right)
=\mbox{ind}\, D_{+,F,\phi(F^\bot),\varepsilon}
=0 ,\eqno(2.44)$$
where $L(F^\bot)$ is the Hirzebruch $L$-class of $F^\bot$.

Now it is a standard fact in topology that any rational Pontrjagin class
 of $F^\bot$ can be expressed as a rational linear combination of the 
classes of the form  
$L(F^\bot)\mbox{ch}\, (\phi(F^\bot))$. So the theorem follows from (2.44). $\Box$

$\ $

{\bf Corollary 2.7.} {\em Under the hyperthesis of Theorem 2.6, one has 
$\langle\widehat{A}(TM), [M]\rangle=0$.}

$\ $

{\bf \S 3. New vanishing theorems} 

\vspace{.2in}

By slightly modifying the construction of the 
operator defined in \S 2, one can also prove the following new vanishing results.

\vspace{ .2in}

{\bf Theorem 3.1}. {\em Let $(M, F)$ be an oriented 
 almost Riemannian foliation with $F$ also oriented. If $M$, 
instead of $F$, is spin, then we still have 
$\langle \widehat{A}(F)p(TM/F),[M]\rangle =0$, under the same condition
that $F$ admits a metric of positive scalar curvature over $M$.  }

\vspace{ .2in}

{\bf Theorem 3.2}. {\em Under the same assumptions as in Theorem 2.6, we have \\
$\langle \widehat{A}(F)e(TM/F),[M]\rangle =0$, 
where $e(TM/F)$ is the Euler class of $TM/F$.}

\vspace{ .2in}

{\bf Appendix. Remarks on the general case} 

$\ $

We first show that the almost isometric foliations studied in [Co]
are almost Riemannian in the sense of \S 2.

Let $(M,F)$ be a foliation, let $TM/F$ be the transversal bundle.
 Let $G$ be the holonomy groupoid of $(M,F)$ (cf. [Co], [W]), then $G$ acts 
on $TM/F$. Let us assume $E$ is a proper subbundle of $TM/F$. We choose 
a splitting of $TM/F$ as 
$$TM/F=E\oplus (TM/F)/E.\eqno(A.1)$$ 
Let $q_1,\, q_2$ be the dimensions of $E$ and $(TM/F)/E$ respectively.

$\ $

{\bf Definition A.1.} If there is a metric  $g^{TM/F}$ on $TM/F$ with 
its restriction to $E$ and $(TM/F)/E$ such that the action of $G$ on $TM/F$ 
takes the form 
$$ \left(\begin{array}{cc}O(q_1)&A\\0&O(q_2)\end{array}\right),
\eqno(A.2)$$
where $O(q_1),\ O(q_2)$ are orthogonal matrices 
of rank $q_1$ and $q_2$ respectively, and $A$ is a $q_1\times q_2$ matrix, 
then we say that $(M,F)$ carries an almost isometric structure.

$\ $

Clearly the existence of the almost isometric structure does not depend on 
the splitting (A.1). Let $g^{TF}$ be a metric on $F$. Choose a splitting 
$$TM=F\oplus F^\bot.\eqno(A.3)$$

We can and we will identify $TM/F$ with $F^\bot$. Thus $E$ and $(TM/F)/E$ are 
identified with subbundles of $F^\bot$ as $F_1^\bot$ and  $F_2^\bot$ 
respectively. Let $g^{F^\bot}$ be the metric on $F^\bot$ corresponding 
to the metric on $g^{TM/F}$, let $g^{F_1^\bot}$ and $g^{F_2^\bot}$ be the 
restrictions of $g^{F^\bot}$ to $F_1^\bot$ and  $F_2^\bot$ respectively. 
Then we have the orthogonal splitting
$$TM=F\oplus F_1^\bot \oplus F_2^\bot,$$
$$g^{TM}=g^F\oplus g^{F_1^\bot}\oplus g^{F_2^\bot}.\eqno(A.4)$$

Recall that the action of the holonomy groupoid $G$ can be defined with 
respect to the canonical flat connection on $F$, the Bott connection 
defined in (1.9). Then the almost isometric condition can be rewritten
 locally as 
$$\langle [X, U_i], V_i\rangle +\langle U_i, [X, V_i]\rangle=X\langle U_i, V_i\rangle,$$
$$\langle [X, U_1], U_2\rangle=0,\eqno(A.5)$$ 
where $X\in 
\Gamma(F),\ U_i, \ V_i\in \Gamma(F_i^\bot)$, $i=1,\ 2$.

Equation (A.5) can be rewritten as 
$$\left\langle X, \nabla_{U_i}V_i+\nabla_{V_i}U_i\right\rangle =0,$$
$$\left\langle \nabla_X U_1,U_2\right\rangle+
\left\langle X,\nabla_{U_1}U_2\right\rangle=0.\eqno(A.6)$$

Thus for any $X\in \Gamma(F)$, $U=U_1+U_2\in \Gamma(F_1^\bot)\oplus 
\Gamma(F_2^\bot)=\Gamma(F^\bot)$, one has, in view of (1.13), that 
$$\omega(X)(U, U)=-{1\over 2 } \left\langle X,\nabla^{TM}_UU\right\rangle
=-{1\over 2 } \left\langle X,\nabla^{TM}_{U_1}U_2+\nabla^{TM}_{U_2}U_1\right\rangle$$
$$=-{1\over 2 } \left\langle U_1,\nabla^{TM}_{X}U_2\right\rangle+{1\over 2 } 
\left\langle U_1,\nabla^{TM}_{U_2}X\right\rangle 
=-{1\over 2 } \left\langle [X,U_2], U_1\right\rangle.\eqno(A.7)$$

Now for any $\gamma>0$, set 
$$g_\gamma^{TM}=g^F\oplus g^{F_1^\bot}\oplus 
{1\over\gamma} g^{F_2^\bot}.$$
Clearly, if $U=U_1+U_2\in \Gamma(F^\bot)$ is of norm one in $g^{TM}$, then  
$U(\gamma)=U_1+\sqrt{\gamma}U_2 \in \Gamma(F^\bot)$ also has norm one in 
$g_\gamma^{TM}$. Let $\omega_\gamma$ be the form as constructed in \S 1,
 corresponding to $g_\gamma^{TM}$. Then by (A.7), we have 
$$\omega_\gamma(X)(U(\gamma),U(\gamma))=-{1\over 2}
\sqrt{\gamma}\left\langle [X,U_2], U_1\right\rangle
\eqno(A.8)$$ 
from which one sees that
$$||\omega_\gamma||_{g_\gamma^{TM}}=\sqrt{\gamma}||\omega||_{g^{TM}}.
\eqno(A.9)$$

Taking $\gamma$ to be as small as possible, we arrive at the following 

$\ $

{\bf Proposition A.2.} {\em  Any almost isometric foliation carries an
 almost Riemannian structure.} 

$\ $

Now recall that in [Co], Connes first proved his vanishing theorem for compact
almost isometric foliations by using the techniques of cyclic cohomology, and then pass
to non-compact manifolds to prove the general case. Thus, what we have done in \S 2 may 
be  thought of as direct geometric approach of the first step of Connes' proof.
As for the general case, it seems one needs a full geometric understanding
of the Connes fibration constructed in [Co], which is  non-compact. 
We leave this for  further studies.

$\ $

{\bf References}

$\ $

{\small [AS] M. F. Atiyah and I. M. Singer, The index of elliptic operators I. 
{\em Ann. of Math.} 87 (1968), 484-530.

[BZ] J.-M. Bismut and W. Zhang, An extension of a theorem by Cheeger and 
M\"{u}ller. {\em Ast\'erisque,} Tom. 205, (1992) Paris.

[Bo] R. Bott, On a topological obstruction to integrability.
 {\em Global Analysis: Proc. Symp. Pure Math.} vol.16, (1970), 127-131.

[Co] A. Connes, Cyclic cohomology and the transverse fundamental class of 
a foliation. {\em Geometric Methods in Operator Algebras.}  H. Araki eds., 
pp. 52-144,
Pitman Res. Notes in Math. Series, vol. 123, 1986.

[CoS] A. Connes and G. Skandalis, The longitudinal index theorem for foliations.
{\it Publ. Res. Inst. Math. Sci. Kyoto} 20 (1984), 1139-1183.

[L] A. Lichnerowicz, Spineurs harmoniques. {\em  C. R. Acad. Sci. Paris, 
S\'erie A,} 257 (1963), 7-9.

[W] E. Winkelnkemper, The graph of a foliation. {\it Ann. Global Anal. 
Geom.} 1 (1983), 51-75.

[Y] S. T. Yau, Private communications. 1992.}

$\ $

K. L.: Department of Mathematics,  
Stanford University, Palo Alto, CA 94305, U. S. A. 

E-mail: kefeng@math.stanford.edu

$\ $

W. Z.: 
Nankai Institute of Mathematics, Tianjin 300071, P. R. China.

E-mail: weiping@nankai.edu.cn

\end{document}